\documentclass{article}

\title{A self-contained, brief and complete formulation of Voevodsky's Univalence Axiom}
\author{Mart{\'\i}n H\"otzel Escard\'o}

\usepackage{amsmath,amssymb,url,hyperref}

\newcommand{\N}{\mathbb{N}}
\newcommand{\Id}{\operatorname{Id}}
\newcommand{\J}{\operatorname{J}}
\newcommand{\refl}{\operatorname{refl}}
\newcommand{\K}{\operatorname{K}}
\newcommand{\U}{\operatorname{\mathcal{U}}}
\newcommand{\V}{\operatorname{\mathcal{V}}}
\newcommand{\isSet}{\operatorname{isSet}}
\newcommand{\Grp}{\operatorname{Grp}}
\newcommand{\isSingleton}{\operatorname{isSingleton}}
\newcommand{\isEquiv}{\operatorname{isEquiv}}
\newcommand{\Eq}{\operatorname{Eq}}
\newcommand{\singletonType}{\operatorname{singletonType}}
\newcommand{\id}{\operatorname{id}}
\newcommand{\idIsEquiv}{\operatorname{idIsEquiv}}
\newcommand{\IdToEq}{\operatorname{IdToEq}}
\newcommand{\isUnivalent}{\operatorname{isUnivalent}}
\newcommand{\Iso}{\operatorname{Iso}}

\begin{document}

\maketitle

\begin{abstract}
  In introductions to the subject for a general audience of
  mathematicians or logicians, the univalence axiom is typically
  explained by handwaving. This gives rise to several misconceptions,
  which cannot be properly addressed in the absence of a precise
  definition. In this short set of notes we give a complete
  formulation of the univalence axiom from scratch. The underlying
  idea of these notes is that they should be as concise as possible
  (and not more). They are not meant to be an Encyclopedia of
  Univalence.

\end{abstract}

\section{Introduction}

The univalence axiom~\cite{unimath,hottbook,grayson} is not true or false
in, say, ZFC or the internal language of an elementary topos. It
cannot even be formulated. As the saying goes, it is not even wrong.
This is because
\begin{quote}
   univalence is a property of Martin-L\"of's \emph{identity type}
   of a universe of types.
\end{quote}
Nothing like Martin-L\"of's identity type occurs in ZFC or topos logic
as a \emph{native} concept. Of course, we can create \emph{models} of
the identity type in these theories, which will make univalence hold
or fail. But in these notes we try to understand the primitive concept
of identity type, directly and independently of any such particular
model, as in (intensional) Martin-L\"of type
theory~\cite{MR0387009,MR1686864}, and the univalence axiom for it. In
particular, we don't use the equality sign ``$=$'' to denote the
identity type $\Id$, or think of it as a path space.

Univalence is a type, and the univalence axiom says that this type has
some inhabitant. It takes a number of steps to construct this type, in
addition to subtle decisions (e.g.\ to work with equivalences rather
than isomorphisms, as discussed below).

We first need to briefly introduce Martin-L\"of type theory (MLTT). We
will not give a full definition of MLTT. Instead, we will mention
which constructs of MLTT are needed to give a complete definition of
the univalence type. This will be enough to illustrate the important
fact that in order to understand univalence we first need to
understand Martin-L\"of type theory well.

\section{Martin-L\"of type theory, briefly}

\subsection{Types and their elements}

Types are the analogues of sets in ZFC and of objects in topos theory.
Types are constructed together with their elements, and not by
collecting some previously existing elements. When a type is
constructed, we get freshly new elements for it. We write
\[
    x:X
\]
to declare that the element $x$ has type $X$. This is not something
that is true or false, unlike a membership relation $x \in X$ in
ZFC. In other words, $x \in X$ in ZFC is a binary relation, whereas
$x:X$ in type theory simply specifies that $x$ ranges over $X$.
For example, if $\N$ is the type of natural numbers, we may
write
\begin{gather*}
    0 : \N, \\
   (0,0) : \N \times  \N.
\end{gather*}
However, the following statements are nonsensical and syntactically
incorrect, rather than false:
\begin{eqnarray*}
    0 : \N \times  \N & \text{(nonsense)}, \\
   (0,0) : \N & \text{(nonsense)}.
\end{eqnarray*}
This is no different from the situation in the internal language of a
topos.

\subsection{Products and sums of type families}

Given a family of types $A(x)$ indexed by elements $x$ of a type $X$, we can
form its product and sum:
\begin{gather*}
    \Pi(x:X), A(x), \\
    \Sigma(x:X), A(x),
\end{gather*}
which we also write $\Pi A$ and $\Sigma A$. An element of the type $\Pi A$ is a
function that maps elements $x:X$ to elements of $A(x)$. An element of the
type $\Sigma A$ is a pair $(x,a)$ with $x:X$ and $a:A(x)$.
(We adopt the convention that $\Pi$ and $\Sigma$ scope over the whole rest of the
expression.)

We also have the type $X\to Y$ of functions from $X$ to $Y$, which is
the particular case of $\Pi$ with the constant family $A(x):=Y$.  The
cartesian product $X\times Y$, whose elements are pairs, is the
particular case of $\Sigma$ with $A(x):=Y$ again.
We also have the disjoint sum $X+Y$, the empty type and the
one-element type, which will not be needed to formulate univalence.

\subsection{Quantifiers and logic}

There is no underlying logic in MLTT. Propositions are types, and
$\Pi$ and $\Sigma$ play the role of universal and existential
quantifiers, via the so-called Curry-Howard interpretation of
logic. As for the connectives, implication is given by the
function-type construction~$\to$, conjunction by the binary cartesian
product~$\times$ , disjunction by the binary disjoint sum~$+$, and
negation by the type of functions into the empty type.

When a type is understood as a proposition, its elements correspond to
proofs.  In this case, instead of saying that $A$ has a given element,
it is common practice to say that $A$ holds. Then a type declaration
$x:A$ is read as saying that $x$ is a proof of $A$. But this is just a
linguistic device, which is not reflected in the formalism.

We remark that in univalent mathematics the terminology
\emph{proposition} is reserved for subsingleton types (types whose
elements are all identified). The propositions that arise in the
construction of the univalence type are all subsingletons.

\subsection{The identity type}

Given a type $X$ and elements $x,y:X$, we have the identity type
\[
    \Id_X(x,y),
\]
with the subscript $X$ often elided. The idea is that $\Id(x,y)$ collects
the ways in which $x$ and $y$ are identified.
We have a function
\[
    \refl : \Pi(x:X), \Id(x,x),
\]
which identifies any element with itself. Without univalence, $\refl$ is
the only given way to construct elements of the identity type.

In addition to $\refl$, for any given type family $A(x,y,p)$
indexed by elements $x,y:X$ and $p:\Id(x,y)$ and any given function
\[
    f : \Pi(x:X), A(x,x,\refl(x)),
\]
we have a function
\[
 \J(A,f) :  \Pi(x,y:X), \Pi(p:\Id(x,y)), A(x,y,p)
\]
with $\J(A,f)(x,x,\refl(x))$ stipulated to be $f(x)$.
We will see examples of uses of $\J$ in the steps leading to the
construction of the univalence type.

Then, in summary, the identity type is given by the data $\Id,\refl,\J$.
With this, the exact nature of the type $\Id(x,y)$ is fairly
under-specified. It is consistent that it is always a subsingleton in
the sense that $\K(X)$ holds, where
\[
   \K(X) := \Pi(x,y:X), \Pi(p,q:\Id(x,y)), \Id(p,q).
\]
The second identity type $\Id(p,q)$ is that of the type
$\Id(x,y)$. This is possible because any type has an identity type,
including the identity type itself, and the identity type of the
identity type, and so on, which is the basis for univalent mathematics
(but this is not discussed here, as it is not needed in order to
construct the univalence type).

The $\K$ axiom says that $\K(X)$ holds for every type $X$.  In
univalent mathematics, a type $X$ that satisfies $\K(X)$ is called a
set, and with this terminology, the $\K$ axiom says that all types are
sets.

On the other hand, the univalence axiom provides a means of
constructing elements other than $\refl(x)$, at least for some types,
and hence the univalence axiom implies that some types are not
sets. (Then they will instead be 1-groupoids, or 2-groupoids, \dots, or
even $\infty$-groupoids, with such notions defined within MLTT rather
than via models, but we will not address this important aspect of
univalent mathematics here).

\subsection{Universes}

Our final ingredient is a ``large'' type of ``small'' types, called a
universe. It is common to assume a tower of universes $\U_0, \U_1, \U_2, \dots
$ of ``larger and larger'' types, with
\begin{gather*}
   \U_0 : \U_1, \\
   \U_1 : \U_2, \\
   \U_2 : \U_3, \\
\vdots
\end{gather*}

When we have universes, a type family $A$ indexed by a type $X:\U$ may be
considered to be a function $A:X\to \V$ for some universe $\V$.
Universes are also used to construct types of mathematical structures,
such as the type of groups, whose definition starts like this:
\[
 \Grp := \Sigma(G:\U), \isSet(G) \times  \Sigma(e:G), \Sigma(-\cdot- : G\times G\to G), (\Pi(x:G), \Id(e \cdot x,x)) \times  \cdots
\]
Here $\isSet(G):=\Pi(x,y:G),\Pi(p,q:\Id(x,y)),\Id(p,q)$, as above. With
univalence, $\Grp$ itself will not be a set, but a 1-groupoid instead,
namely a type whose identity types are all sets. Moreover, if $\U$
satisfies the univalence axiom, then for $A,B:\Grp$, the identity type
$\Id(A,B)$ can be shown to be in bijection with the group isomorphisms of
$A$ and $B$.

\section{Univalence}

Univalence is a property of the identity type $\Id_{\U}$ of a universe $\U$. It
takes a number of steps to define the univalence type.

\subsection{Construction of the univalence type}

We say that a type $X$ is a \emph{singleton} if we have an element $c:X$ with
$\Id(c,x)$ for all $x:X$. In Curry-Howard logic, this is
\[
    \isSingleton(X) := \Sigma(c:X), \Pi(x:X), \Id(c,x).
\]
For a function $f:X\to Y$ and an element $y:Y$, its fiber is the type of
points $x:X$ that are mapped to (a point identified with) $y$:
\[
    f^{-1}(y) := \Sigma(x:X),\Id(f(x),y).
\]
The function $f$ is called an equivalence if its fibers are all
singletons:
\[
    \isEquiv(f) := \Pi(y:Y), \isSingleton(f^{-1}(y)).
\]
The type of equivalences from $X:\U$ to $Y:\U$ is
\[
    \Eq(X,Y) := \Sigma(f:X\to Y), \isEquiv(f).
\]
Given $x:X$, we have the singleton type consisting of the elements $y:X$
identified with $x$:
\[
   \singletonType(x) := \Sigma(y:X), \Id(y,x).
\]
We also have the element $\eta(x)$ of this type:
\[
   \eta(x) := (x, \refl(x)).
\]
We now need to \emph{prove} that singleton types are singletons:
\[
   \Pi(x:X), \isSingleton(\singletonType(x)).
\]
In order to do that, we use $\J$ with the type family
\[
   A(y,x,p) := \Id(\eta(x),(y,p)),
\]
and the function
\begin{eqnarray*}
  f & : & \Pi(x:X), A(x,x,\refl(x)) \\
   f(x) & := & \refl(\eta(x)).
\end{eqnarray*}
With this we get a function
\begin{eqnarray*}
   \phi & : & \Pi(y,x:X), \Pi(p:\Id(y,x)), \Id(\eta(x),(y,p)) \\
   \phi & := & \J(A,f).
\end{eqnarray*}
Notice the reversal of $y$ and $x$.
With this, we can in turn define a function
\begin{eqnarray*}
   g & : & \Pi(x:X), \Pi(\sigma :\singletonType(x)), \Id(\eta(x),\sigma ) \\
   g(x,(y,p)) & := & \phi(y,x,p).
\end{eqnarray*}
Finally, using the function $g$, we get our desired result, that
singleton types are singletons:
\begin{eqnarray*}
   h & : & \Pi(x:X), \Sigma(c:\singletonType(x)), \Pi(\sigma :\singletonType(x)), \Id(c,\sigma )\\
   h(x) & := & (\eta(x),g(x)).
\end{eqnarray*}
Now, for any type $X$, its identity function $\id_X$, defined by
$\id(x) := x$, is an equivalence. This is because the fiber
$\id^{-1}(x)$ is simply the singleton type defined above, which we
proved to be a singleton. We need to name this function:
\[
   \idIsEquiv : \Pi(X:\U), \isEquiv(\id_X).
 \]
The identity function $\id_X$ should not be confused with the identity
type $\Id_X$. Now we use $\J$ a second time to define a function
\[
   \IdToEq : \Pi(X,Y:\U), \Id(X,Y) \to  \Eq(X,Y).
\]
For $X,Y:\U$ and $p:\Id(X,Y)$, we set
\[
   A(X,Y,p) := \Eq(X,Y)
\]
and
\[
   f(X) := (\id_X , \idIsEquiv(X))
\]
and
\[
   \IdToEq := \J(A,f).
\]
Finally, we say that the universe $\U$ is univalent if the map
$\IdToEq(X,Y)$ is itself an equivalence:
\[
   \isUnivalent(\U) := \Pi(X,Y:\U), \isEquiv(\IdToEq(X,Y)).
\]

\subsection{The univalence axiom}

The type $\isUnivalent(\U)$ may or may not have an inhabitant. The
univalence axiom says that it does. The $\K$ axiom implies that it doesn't.
Because both univalence and the $\K$ axiom are consistent, it follows that univalence is undecided in MLTT.

\subsection{Notes}

\begin{enumerate}
\item
 The minimal Martin-L\"of type theory needed to formulate univalence
    has
    \[
      \Pi, \Sigma, \Id, \U, \U'.
    \]
    Two universes $\U :\U'$ suffice, where univalence talks about $\U$.

\item It can be shown, by a very complicated and interesting argument,
    that
    \[
     \Pi(u,v: \isUnivalent(\U)), \Id(u,v).
   \]
    This says that univalence is a subsingleton type (any two of its
    elements are identified). In the first step we use $u$ (or $v$) to get
    function extensionality (any two pointwise identified functions
    are identified), which is \emph{not} provable in MLTT, but is provable
    from the assumption that $\U$ is univalent. Then, using this, one
    shows that being an equivalence is a subsingleton type. Finally,
    again using function extensionality, we get that a product of
    subsingletons is a subsingleton. But then $\Id(u,v)$ holds, which is
    what we wanted to show. But this of course omits the proof that
    univalence implies function extensionality (originally due to
    Voevodsky), which is fairly elaborate.

 \item For a function $f:X\to Y$, consider the type
   \[
     \Iso(f) := \Sigma(g:Y\to X), (\Pi(x:X), \Id(g(f(x)),x)) \times  (\Pi(y:Y), \Id(f(g(y)),y)).
   \]
    We have functions $r:\Iso(f)\to \isEquiv(f)$ and
    $s:\isEquiv(f)\to \Iso(f)$. However, the type $\isEquiv(f)$ is always a
    subsingleton, assuming function extensionality, whereas the type
    $\Iso(f)$ need not be. What we do have is that the function $r$ is a
    retraction with section $s$.

    Moreover, the univalence type formulated as above, but using
    $\Iso(f)$ rather than $\isEquiv(f)$, is provably empty, e.g.\ for
    MLTT with $\Pi, \Sigma, \Id$, the empty and two-point types, and
    two universes, as shown by Shulman~\cite{shulman:e46}. With only
    one universe, the formulation with $\Iso(f)$ is consistent, as
    shown by Hofmann and Streicher's groupoid model~\cite{MR1686862},
    but in this case all elements of the universe are sets and
    $\Iso(f)$ is a subsingleton, and hence equivalent to
    $\isEquiv(f)$.

    So, to have a consistent axiom in general, it is crucial to use
    the type $\isEquiv(f)$. It was Voevodsky's insight not only that a
    subsingleton version of $\Iso(f)$ is needed, but also how to
    construct it. The construction of $\isEquiv(f)$ is very simple and
    elegant, and motivated by homotopical models of the theory, where
    it corresponds to the concept with the same name. But the
    univalence axiom can be understood without reference to homotopy
    theory.

  \item Voevodsky gave a model of univalence for MLTT with
    $\Pi,\Sigma$, empty type, one-point type, two-point type, natural
    numbers, and an infinite tower of universes in simplicial
    sets~\cite{kapulkin:lumsdaine:voevodsky,kapulkin:lumsdaine}, thus
    establishing the consistency of the univalence axiom.

    The consistency of the univalence axiom shows that, before we
    postulate it, MLTT is ``proto-univalent'' in the sense that it
    cannot distinguish concrete isomorphic types such as $X:=\N$ and
    $Y:=\N\times \N$ by a property $P:\U\to \U$ such that $P(X)$ holds
    but $P(Y)$ doesn't.  This is because, being isomorphic, $X$ and
    $Y$ are equivalent. But then univalence implies $\Id(X,Y)$, which
    in turn implies $P(X) \iff P(Y)$ using $\J$.  Because univalence
    is consistent, it follows that for any given concrete
    $P:\U\to \U$, it is impossible to prove that $P(X)$ holds but
    $P(Y)$ doesn't.

    So MLTT is invariant under isomorphism in this doubly negative,
    meta-mathematical sense. With univalence, it becomes invariant
    under isomorphism in a positive, mathematical sense.

 \item Thus, we see that the formulation of univalence is far from
    direct, and has much more to it than the (in our opinion,
    misleading) slogan ``isomorphic types are equal''.

    What the consistency of the univalence axiom says is that one
    possible understanding of Martin-L\"of's identity type $\Id(X,Y)$ for
    $X,Y:\U$ is as precisely the type $\Eq(X,Y)$ of equivalences, in
    the sense of being in one-to-one correspondence with it. Without
    univalence, the nature of the identity type of the universe in
    MLTT is fairly under-specified. It is a remarkable property of
    MLTT that it is consistent with this understanding of the identity
    type of the universe, discovered by Vladimir Voevodsky (and
    foreseen by Martin Hofmann and Thomas Streicher~\cite{MR1686862}
    in a particular case).
  \item It should also be emphasized that what univalence does it to
    express the identity type $\Id(X,Y)$ for $X,Y : \U$ in terms of
    the identity types of the types $X$ and $Y$.  This is because the
    notion of equivalence $X \simeq Y$ is defined in terms of the
    identity types of $X$ and $Y$. In this sense, univalence is an
    extensionality axiom: it says what identity of types $X$ and $Y$
    is in terms of what identity for the elements of the types $X$ and
    $Y$ are.  From this perspective, it is very interesting that
    univalence implies function extensionality (any two pointwise
    identified functions are themselves identified) and propositional
    extensionality (any two subsingletons, or truth values, which
    imply each other are identified). Thus, univalence is a common
    generalization of function extensionality and propositional
    extensionality.
  \end{enumerate}
  This paper only explains what the \emph{univalence axiom} is. A
  brief and reasonably complete introduction to \emph{univalent
    mathematics} is given by Grayson~\cite{grayson}.

  \section*{Acknowledgements}

  I benefitted from input by Andrej Bauer, Marta Bunge, Thierry
  Coquand, Dan Grayson and Mike Shulman on draft versions of these
  notes.

\nocite{escardo:ufs}
\bibliographystyle{plain}
\bibliography{UnivalenceFromScratch}

\begin{thebibliography}{10}

\bibitem{escardo:ufs}
Mart\'{\i}n~H\"otzel Escard\'o.
\newblock Univalence from scratch in {A}gda, February 2018.
\newblock
  \url{https://arxiv.org/src/1803.02294v1/anc/UnivalenceFromScratch.lagda}.

\bibitem{grayson}
Daniel~R. Grayson.
\newblock An introduction to univalent foundations for mathematicians.
\newblock {\em Bull. Amer. Math. Soc.}, 2018.
\newblock \url{https://doi.org/10.1090/bull/1616}.

\bibitem{MR1686862}
Martin Hofmann and Thomas Streicher.
\newblock The groupoid interpretation of type theory.
\newblock In {\em Twenty-five years of constructive type theory ({V}enice,
  1995)}, volume~36 of {\em Oxford Logic Guides}, pages 83--111. Oxford Univ.
  Press, New York, 1998.

\bibitem{kapulkin:lumsdaine}
Chris Kapulkin and Peter~LeFanu Lumsdaine.
\newblock The simplicial model of univalent foundations (after {V}oevodsky),
  2012.
\newblock arXiv:1211.2851.

\bibitem{kapulkin:lumsdaine:voevodsky}
Chris Kapulkin, Peter~LeFanu Lumsdaine, and Vladimir Voevodsky.
\newblock The simplicial model of univalent foundations, 2012.
\newblock arXiv:1203.2553.

\bibitem{MR0387009}
Per Martin-L\"of.
\newblock An intuitionistic theory of types: predicative part.
\newblock pages 73--118. Studies in Logic and the Foundations of Mathematics,
  Vol. 80, 1975.

\bibitem{MR1686864}
Per Martin-L\"of.
\newblock An intuitionistic theory of types.
\newblock In {\em Twenty-five years of constructive type theory ({V}enice,
  1995)}, volume~36 of {\em Oxford Logic Guides}, pages 127--172. Oxford Univ.
  Press, New York, 1998.

\bibitem{shulman:e46}
Michael Shulman.
\newblock Solution to {E}xercise 4.6 (in pure {MLTT}), May 2018.
\newblock
  \newline\url{https://github.com/HoTT/HoTT/blob/master/contrib/HoTTBookExercises.v}.

\bibitem{hottbook}
The {Univalent Foundations Program}.
\newblock {\em Homotopy Type Theory: Univalent Foundations of Mathematics}.
\newblock \url{https://homotopytypetheory.org/book}, Institute for Advanced
  Study, 2013.

\bibitem{unimath}
Vladimir Voevodsky.
\newblock An experimental library of formalized mathematics based on the
  univalent foundations.
\newblock {\em Math. Structures Comput. Sci.}, 25(5):1278--1294, 2015.

\end{thebibliography}

\end{document}